\documentclass[letterpaper, 10 pt, conference]{ieeeconf}  

\IEEEoverridecommandlockouts                              
\usepackage{graphicx} 
\usepackage{amsmath} 
\usepackage{amssymb} 
\usepackage{cite}
\usepackage{algorithm}
\usepackage{algorithmic}
\usepackage{subfigure}
\usepackage{hyperref}

\newcommand{\diff }[2]{\frac{    d   {#1}}{    d   {#2}}}
\newcommand{\pdiff}[2]{\frac{\partial{#1}}{\partial{#2}}}

\title{\LARGE \bf
Filtering and Prediction of the Blood Glucose Concentration\\ using an Android Smart Phone and a Continuous Glucose Monitor
}

\author{Zeinab Mahmoudi, Dimitri Boiroux, Tobias K. S. Ritschel, John Bagterp J{\o}rgensen
	\thanks{*Z. Mahmoudi, D. Boiroux, Tobias K. S. Ritschel and  J. B. J{\o}rgensen are with the Department of Applied Mathematics and Computer Science, Technical University of Denmark, DK-2800 Kgs. Lyngby, Denmark. Corresponding author: J. B. J{\o}rgensen (E-mail: {\tt\small jbjo@dtu.dk})}
}

\begin{document}
\maketitle
\thispagestyle{empty}
\pagestyle{empty}

\begin{abstract}
In this paper we numerically assess the performance of Java linear algebra libraries for the implementation of nonlinear filters in an Android smart phone (Samsung A5 2017). We implemented a linear Kalman filter (KF), an extended Kalman filter (EKF), and an unscented Kalman filter (UKF). These filters are used for state and parameter estimation, as well as fault detection and meal detection in an artificial pancreas. We present the state estimation technologies used for glucose estimation based on a continuous glucose monitor (CGM). We compared three linear algebra libraries: The Efficient Java Matrix Library (EJML), JAMA and Apache Common Math. Overall, EJML provides the best performance for linear algebra operations. We demonstrate the implementation and performance of filtering (KF, EKF and UKF) using real CGM data.
\end{abstract}

\section{Introduction}
\label{sec:Introduction}

Kalman filtering is used to estimate the states of a linear state-space model where discrete-time measurements are available. Variations of the Kalman filter also exist for nonlinear systems. The extended Kalman filter (EKF) uses a first order approximation to update the mean and covariance of the state vector. In the unscented Kalman filter (UKF) \cite{Julier:etal:2000,Sarkka:2007}, we select a number of representative points around the mean value (also called sigma points) to numerically evaluate the mean and covariance of the state vector. The particle filter (PF) and ensemble Kalman filter (EnKF) propagate a large number of points to estimate the probability distribution of the states \cite{Burgers:etal:1998,Tulsyan:etal:2016}. Kolmogorov forward equations (also called Fokker-Planck equations) provide an analytical expression of the state probability density function \cite{Challa:etal:2000}. 

When the measurement sampling time is fast enough compared to the nonlinearities, the EKF and the UKF provide a suitable approximation of the true system. The PF and the EnKF are more accurate than the EKF and the UKF and can capture non-Gaussian distributions, but suffer from the curse of dimensionality. The Kolmogorov forward equations are exact but need to solve a system of partial differential equations and cannot be numerically solved in a reasonable time when the number of states become too large. The artificial pancreas (AP) is an example of a closed-loop control system. Filtering is important to ensure good performance of the closed-loop algorithm \cite{Cobelli:etal:2011,Trevitt:etal:2016,Boiroux:etal:NOLCOS2010,Boiroux:etal:2010,Schmidt:etal:2013,Boiroux:etal:BSPC2018,Boiroux:etal:2017:CEP,Boiroux:etal:BMS2018NMPC,Boiroux:etal:CDC2018}. The AP gets measurements from a continuous glucose monitor (CGM), computes the optimal insulin dose, and sends this dose to a continuous subcutaneous insulin infusion (CSII) pump. In the AP, filtering and prediction are used for model identification \cite{Boiroux:etal:DYCOPS2016}, detection of meals \cite{Mahmoudi:etal:2019}, detection of faulty CGM measurements \cite{mahmoudi2017a}, as well as in the model predictive control (MPC) \cite{Boiroux:etal:BMS2018NMPC} algorithm. The filtering is also important for CGM enabled insulin pen systems \cite{Boiroux:etal:2015:BMS,Boiroux:etal:JDST2017}. For the AP, EKF and UKF can be used as filtering algorithms and perform better than the PF \cite{Mahmoudi:etal:EMBC2016} due to the fast CGM sampling time, typically every five minutes. More generally speaking, filtering, detection and numerical optimization in the AP involve linear algebra operations for small matrices. In this paper, we evaluate the performance of linear KF and nonlinear filters (EKF and UKF) implemented in an Android platform. We compare three Java linear algebra libraries and demonstrate the implementation of KF, EKF and UKF.

The rest of the paper is structured as follows. Section \ref{sec:LinearModel} presents the linear Kalman filter (KF), and Section \ref{sec:NonlinearModel} explains the formulation for the extended Kalman filter (EKF), and unscented Kalman filter (UKF). Section \ref{sec:Implementation} describes the Java implementation of the filter in an Android smart phone, and the hardware configuration. The results are presented in Section \ref{sec:Results}. Section \ref{sec:Conclusions} summarizes the main contributions of this paper.
\section{Linear Model}\label{sec:LinearModel}
We use is the Medtronic virtual patient (MVP) type 1 diabetes patient model described in \cite{kanderian2009a}. The MVP is a nonlinear model and we use it for the nonlinear KF. For the linear KF, we identify a  second-order transfer function model for insulin to sc glucose dynamics and for the CHO to sc glucose dynamics. The model in the Laplace domain is defined as
\begin{subequations}
\begin{equation}
Z(s)=G_u(s)U(s)+G_d(s)D(s),
\end{equation}
where 
\begin{align}
&G_u(s)=\dfrac{K_u}{(\tau_us+1)^2}, \label{ins}\\
&G_d(s)=\dfrac{K_d}{(\tau_ds+1)^2}. \label{d}
\end{align}
\label{eq:eqmodel}
\end{subequations}
The output, $Z(s)$, is the sc glucose concentration as the deviation variable, which is the deviation from the steady-state sc glucose concentration (100 mg/dL). $U(s)$ is the sc insulin input rate (IU/min), and $D(s)$ is the CHO ingestion rate (g/min). $U$ and $D$ are also deviation variables. The transfer functions, $G_u(s)$ and $G_d(s)$, are the Laplace transforms of the insulin and CHO impulse responses, respectively. The gains, $K_I$ and $K_d$, correspond to the steady state change in BG for a unit step in the inputs, and the time constants, $\tau_u$ and $\tau_d$, determine the time to reach the steady state. Mahmoudi et al \cite{Mahmoudi:etal:2019} present a method for identification of the parameters, $K_u$, $K_d$, $\tau_u$, and $\tau_d$. The identified model is then converted to a linear time-invariant discretized state-space model \cite{Mahmoudi:etal:2019}, which is 
\begin{subequations}
	\label{eq:LinearSystemXZY}
\begin{alignat}{3}
	x_{k+1} &= A x_k + B u_k + E (d_k + w_k), \\
	z_k &= C x_k, \\
	y_k &= z_k + v_k,
\end{alignat}
\label{lm}
\end{subequations}
with 
\begin{align}
	w_k &\sim N_{iid}(0,Q), \\ 
	v_k &\sim N_{iid}(0,R),
\end{align}
where $u_k$ is the subcutaneous insulin infusion rate, $d_k$ is the meal ingestion, $y_k$ is the CGM sensor measurement, and $v_k$ is the sensor noise \cite{biagi2017a}. For the linear model, We assume that the uncertainty enters the process through the meal ingestion and therefore, we add the process noise, $w$, to the CHO input and use $Q$ as a tuning parameter \cite{Mahmoudi:etal:ccta}. For the linear KF implementation, we compute the stationary state covariance, $P_{\infty|\infty}$ and the stationary filter gain as follows.

The discrete algebraic Riccati equation (DARE),
\begin{equation}
\begin{split}
	P &= A  P A^T + EQE^T \\ & \qquad - (A P C^T ) (C P C^T + R)^{-1} (A P C^T)^T,
\end{split}
\end{equation}
for the KF of the discrete-time stochastic linear model \eqref{lm} may be used to compute the stationary covariance matrix of the one-step prediction, $P = \lim_{k \rightarrow \infty} P_{k|k-1} $. The corresponding matrices $R_{e,\infty} = \lim_{k \rightarrow \infty} R_{e,k}$, $K_\infty = \lim_{k \rightarrow \infty} K_k$, and $P_{\infty|\infty} = \lim_{k \rightarrow \infty} P_{k|k}$ are computed by
\begin{subequations}
\begin{align}
	R_{e,\infty} &= C P C^T + R, \\
	K_{\infty} &= P C^T R_{e,\infty}^{-1}, \\
	P_{\infty|\infty} &= P - K_{\infty} R_{e,\infty} K_{\infty}^T.
\end{align}
$R_{e,\infty}$ is the innovation covariance.
	\begin{align}
	\label{eq:kf:gain:cov}
	R_{e,k} &= C P_{k|k-1} C' + R.
	\end{align}
\end{subequations}  

\subsection{Kalman Filter}
The KF is initialized with $\hat{x}_{0|-1} = x_0$.  
\subsubsection{One-step prediction}
The one-step predictions of the states and the covariance are
\begin{subequations}\label{predkf}
	\begin{align}
	\hat{x}_{k|k-1} &= A \hat{x}_{k-1|k-1} + B u_{k-1} + E d_{k-1},
	\end{align}
\end{subequations}
and the one-step predictions of the outputs and the measurements are
\begin{subequations}
	\begin{align}
	\hat{z}_{k|k-1} &= C \hat{x}_{k|k-1}, \\
	\hat{y}_{k|k-1} &= \hat{z}_{k|k-1}.
	\end{align}
\end{subequations}

\subsubsection{Measurement update}
When the new measurements, $y_k$, become available, we compute the innovation,
\begin{align}
	e_k = y_k - \hat y_{k|k-1}, 
\end{align}
as well as the filtered states and the corresponding covariance matrix,
%
	\begin{align}
	\label{eq:kf:filter:state}
	\hat x_{k|k} 	&= \hat x_{k|k-1} + K_{\infty} e_k.
	\end{align}

\section{Nonlinear Model}\label{sec:NonlinearModel}
In this section, we describe the EKF and the UKF for stochastic differential equations (SDEs) in the form
\begin{align}\label{eq:sde}
	dx(t)
	&= f(x(t), u(t), d(t), p) dt  
	 + g(x(t), u(t), d(t), p) d\omega(t).
\end{align}
The model, $f$, is the MVP model \cite{kanderian2009a}, and in our application the function, $g(x(t), u(t), d(t), p)$, is an additive constant diffusion coefficient and therefore the model reduces to 
\begin{align}\label{eq:sde}
	dx(t)
	&= f(x(t), u(t), d(t), p) dt  + \sigma d\omega(t).
\end{align}
As in Section \ref{sec:LinearModel}, $x(t)$ are the states, $u(t)$ are the manipulated inputs, $d(t)$ are disturbance variables, and $p$ are parameters. $\omega(t)$ is a standard Wiener process, i.e. the incremental covariance of $\omega(t)$ is $Idt$. At the initial time, the states are normally distributed: $x(t_0)\sim N(x_0, P_0)$. The first term in \eqref{eq:sde} is the drift term which represents the deterministic part of the model equations. The second term is the diffusion term which represents the uncertainty in the process, i.e. the process noise. 
 The process outputs,
\begin{equation}
	z(t) = h(x(t), p),
\end{equation}
are measured at discrete points in time, $t_k$:
\begin{equation}
	y(t_k) = z(t_k) + v(t_k).
\end{equation}
$v(t_k)\sim N(0, R)$ is the measurement noise.

\subsection{Extended Kalman Filter}
The EKF is initialized using the distribution of the initial states, i.e. $\hat{x}_{0|-1} = x_0$ and $P_{0|-1} = P_0$.

\subsubsection{One-step prediction}
The one-step prediction of the states and their covariance are the solutions to
\begin{subequations}
	\begin{align}
	\hat x_{k-1}(t_{k-1}) 		&= \hat x_{k-1|k-1}, \\
	P_{k-1}(t_{k-1}) 			&= P_{k-1|k-1}, \\
	\diff{}{t} \hat x_{k-1}(t) 	&= f(\hat{x}_{k-1}(t), u(t), d(t), p), \\
	\diff{}{t} P_{k-1}(t) 		&= A_k(t) P_{k-1}(t) + P_{k-1}(t) A_k(t)' + Q_k(t),
	\end{align}
\end{subequations}
for $t\in[t_{k-1}, t_k]$ where
\begin{subequations}
	\begin{align}
	A_k(t) &= \pdiff{f}{x}(\hat x_{k-1}(t), u(t), d(t), p), \\
	Q_k(t) &=  \sigma  \sigma'.
	\end{align}
\end{subequations}
The one-step predictions are
\begin{subequations}
	\begin{align}
	\hat x_{k|k-1} 	&= \hat x_{k-1}(t_k), \\
	P_{k|k-1} 		&= P_{k-1}(t_k),
	\end{align}
\end{subequations}
and the one-step prediction of the measurements is
\begin{align}
\hat y_{k|k-1} &= \hat z_{k|k-1} = h(\hat x_{k|k-1}, p).
\end{align}
The covariance of the innovation and the Kalman gain are
\begin{subequations}
	\begin{align}
	R_{e, k} 	&= C_k P_{k|k-1} C_k' + R, \\
	K_{fx, k} 	&= P_{k|k-1} C_k' R_{e, k}^{-1},
	\end{align}
\end{subequations}
where $C_k$ is the Jacobian of the sensor model:
\begin{align}
C_k &= \pdiff{h}{x}(\hat{x}_{k|k-1}, p).
\end{align}

\subsubsection{Measurement update}
When the measurements, $y_k$, become available, we compute the innovation,
\begin{align}
	e_k &= y_k - \hat y_{k|k-1},
\end{align}
and the filtered estimates of the states and its covariance:
\begin{subequations}
	\begin{align}
	\hat x_{k|k}
	&= \hat x_{k|k-1} + K_{fx,k} e_k, \\
	P_{k|k}
	&= (I - K_{fx,k} C_k ) P_{k|k-1} (I - K_{fx,k} C_k )' \nonumber \\ 
	&+ K_{fx,k} R K_{fx,k}'.
	\end{align}
\end{subequations}
We use the Joseph-stabilized form of the covariance update in order to ensure positive definiteness of the updated covariance matrix as well as numerical stability.

\subsection{Unscented Kalman Filter}
As for the EKF, the UKF is initialized with $\hat{x}_{0|-1} = x_0$ and $P_{0|-1} = P_0$. 
\subsubsection{One- step prediction}

First, we introduce the weights as \cite{Julier:2002}
\begin{subequations}
\begin{align}
&\lambda=\alpha^2(n+k)-n,  \quad \quad \quad \quad \quad \quad  c=\alpha^2(n+k), \\
& W_m^{(0)}=\lambda/(n+\lambda),\\
& W_c^{(0)}=\lambda/(n+\lambda)+(1-\alpha^2+\beta),\\
& W_m^{(i)}=1/\{2(n+\lambda)\}, \quad \quad \quad \quad \quad i=1,\dots, 2n,\\
& W_c^{(i)}=1/\{2(n+\lambda)\}, \quad \quad \quad \quad \quad i=1,\dots, 2n,\\
& W_m=[W_m^{(0)}\dots W_m^{(2n)}]^T, \\
& W_c=[W_c^{(0)}\dots W_c^{(2n)}],
\end{align}
where $n$ is the number of states in the state-space model.  
 We use $\alpha = 0.01$, $k = 0$, and $\beta = 2$ in this work. Next, the sigma points are generated from the previous filtered states and covariance.
\begin{align}
& \hat{x}^{(0)}_{k-1|k-1}= \hat{x}_{k-1|k-1} \\
& \hat{x}^{(i)}_{k-1|k-1}= \hat{x}_{k-1|k-1}+ \sqrt{c} \Big(\sqrt{P_{k-1|k-1}}\Big)_i,  i=1,\dots , n\\
& \hat{x}^{(i+n)}_{k-1|k-1}= \hat{x}_{k-1|k-1} - \sqrt{c} \Big(\sqrt{P_{k-1|k-1}}\Big)_i,  i= 1,\dots , n.
\end{align}
We use a Cholesky factorization to compute $\sqrt{P_{k-1|k-1}}$, and $\Big(\sqrt{P_{k-1|k-1}}\Big)_i$ denotes the $i$'th column of $\sqrt{P_{k-1|k-1}}$ \cite{Julier:etal:2000}. For each of the sigma points we solve the ODE,
\begin{align}
& \frac{d}{dt} \hat x^{(i)}_{k-1}(t) = f(\hat x^{(i)}_{k-1}(t),u(t),d(t),p), \quad i= 0, \dots, 2n, \nonumber \\
& \quad \quad  \quad  \quad \quad \quad \quad  \quad \quad \quad  \quad \quad \quad \quad \quad t\in[t_{k-1}\ t_{k}].
\end{align}
In parallel, the covariance matrix should be also propagated in time according to
\begin{align}
&\frac{d{P}_{k-1}(t)}{dt}= \nonumber \\
& \sum_{i=0}^{2n}W_c^{(i)}\left( \hat{x}_{k-1}^{(i)}(t)-\hat{x}_{k-1}(t)\right)\left( F^{(i)}(t)-F(t)\right) ^T \nonumber\\ 
& +\sum_{i=0}^{2n}W_c^{(i)} \left(F^{(i)}(t)-F(t)\right) \left( (\hat{x}_{k-1}^{(i)}(t)- \hat{x}_{k-1}(t)\right) ^T \nonumber \\
& \quad +\sigma\sigma^{T}, \quad \quad \quad \quad \quad \quad \quad \quad \quad \quad \quad  t\in[t_{k-1}\ t_{k}],
\end{align}
where,
\begin{align}
& \hat{x}_{k-1}(t)=\sum_{i=0}^{2n}W_m^{(i)}\ \hat{x}_{k-1}^{(i)}(t),  \\
& F^{(i)}(t) = f(\hat{x}_{k-1}^{(i)}(t),u(t),d(t),p),\\
& F(t)=\sum_{i=0}^{2n}W_m^{(i)}\ F^{(i)}(t).
\end{align}
The predicted mean and covariance of the states are computed as 
\begin{align}
&\hat{x}^{(i)}_{k|k-1}=\hat{x}^{(i)}_{k-1}(t_{k}), \quad \quad \quad \quad \quad \quad i= 0, \dots, 2n, \\
\label{xi}
& \hat{x}_{k|k-1}= \sum_{i=0}^{2n}W_m^{(i)}\ \hat{x}_{k|k-1}^{(i)},\\
& P_{k|k-1}= P_{k-1}(t_{k}).
\end{align}
\end{subequations}
\begin{subequations}
 
An intermediate set of sigma points, $\{\hat{\tilde{x}}^{(i)}_{k|k-1}\}_{i=0}^{2n}$, can be generated using the predicted state mean and covariance, 
 \begin{align}
& \hat{\tilde{x}}^{(0)}_{k|k-1}= \hat{x}_{k|k-1}, \\
& \hat{\tilde{x}}^{(i)}_{k|k-1}= \hat{x}_{k|k-1}+ \sqrt{c} \Big(\sqrt{P_{k|k-1}}\Big)_i,  i=1,\dots , n\\
& \hat{\tilde{x}}^{(i+n)}_{k|k-1}= \hat{x}_{k|k-1} - \sqrt{c} \Big(\sqrt{P_{k|k-1}}\Big)_i,  i= 1,\dots , n.
\end{align}
The predicted outputs corresponding to the sigma points are computed according to
\begin{align}
& \hat y^{(i)}_{k|k-1} = \hat z^{(i)}_{k|k-1} = h(\hat{\tilde{x}}^{(i)}_{k|k-1},p),\quad \quad i= 0, \dots, 2n,
\label{yi}
\end{align}
and the mean of the predicted output is
\begin{align}
&\hat{y}_{k|k-1}= \sum_{i=0}^{2n}W_m^{(i)} \hat{y}^{(i)}_{k|k-1}. 
\end{align}
Generation of the intermediate set of sigma points can be omitted for the sake of lowering the computational effort \cite{simon2006a}. In the UKF implementation, we omitted this step and instead we used the sigma points, $\hat{x}^{(i)}_{k|k-1}$, in \eqref{xi} to compute $\hat y^{(i)}_{k|k-1}$ in \eqref{yi}.

The innovation and its covariance are defined as
\begin{align}
& R_{e,k}=\sum_{i=0}^{2n}W_c^{(i)}\left(\hat y^{(i)}_{k|k-1}-\hat{y}_{k|k-1}\right) \left(\hat y^{(i)}_{k|k-1}-\hat{y}_{k|k-1}\right)^T \nonumber \\
& \quad \quad \quad +R. 
\end{align}
The cross-covariance between the predicted states and outputs is
\begin{align}
&R_{xy,k}=\sum_{i=0}^{2n}W_c^{(i)}\ \left( {\hat{x}}^{(i)}_{k|k-1}-\hat{x}_{k|k-1}\right) \left( \hat{y}^{(i)}_{k|k-1}-\hat{y}_{k|k-1}\right)^T.
\label{rxy}
\end{align}
The filter gain is calculated according to 
\begin{align}
&K_{fx,k}=R_{xy,k}\ R_{e,k}^{-1},
\label{k}
\end{align}

\subsubsection{ Measurement update}
 When the measurement $y_{k}$ is available, we compute the innovation,
\begin{align}
&e_{k}=y_{k}-\hat{y}_{k|k-1},
\end{align}
and the filtered state mean and covariance are
\begin{align}
\hat{x}_{k|k} &=\hat{x}_{k|k-1}+K_{fx,k} e_k, \\
\begin{split}
	P_{k|k} &= (I - K_{fx,k} C_k ) P_{k|k-1} (I - K_{fx,k} C_k )' 
	\\ & \qquad + K_{fx,k} R K_{fx,k}'.
\end{split}
\end{align}
\end{subequations}

As in the EKF, we use the Joseph-stabilized form of the covariance update.

\section{Implementation}
\label{sec:Implementation}

\subsection{Java implementation }
We test three linear algebra libraries: JAMA 1.0.3 (\url{http://math.nist.gov/javanumerics/jama/}), EJML 0.36 (\url{http://ejml.org}), and Apache Commons Math 3.6.1 (\url{http://commons.apache.org/proper/commons-math/}). We evaluate the performance of these libraries by comparison of the computation times for the main operations required in the KF, EKF and UKF: Matrix-matrix multiplications, Cholesky decomposition, matrix scaling, and matrix-matrix addition. The prediction and measurement update in the filters require matrix-matrix multiplication (all filters), Cholesky decomposition (EKF and UKF), matrix scaling (KF and UKF), and matrix-matrix addition (all filters). 

\begin{figure}[tb]
\begin{center}
\includegraphics[width=1.03\linewidth]{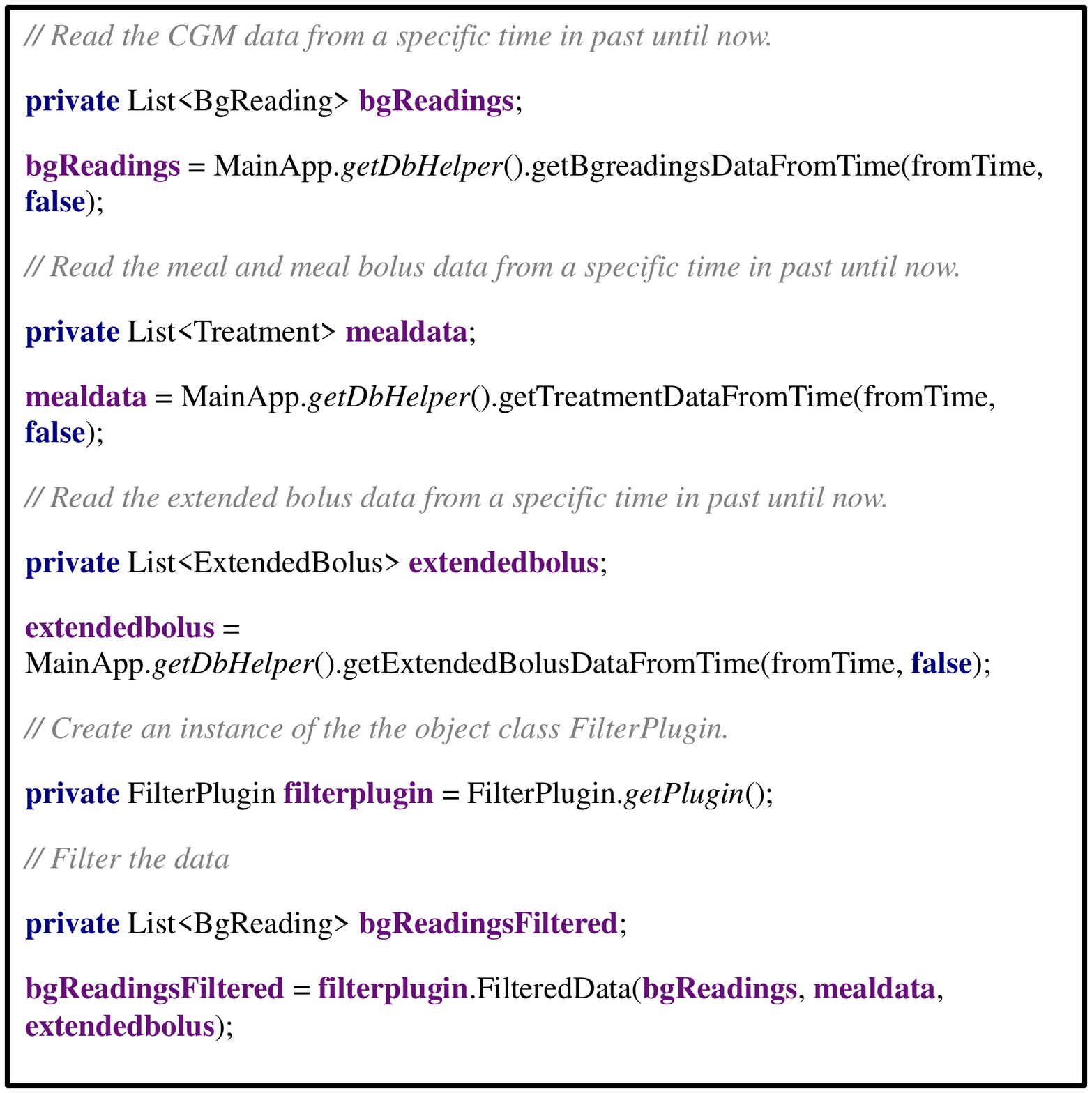}
\vspace{-4cm}
\caption{Java implementation for reading the inputs and filtering the CGM data. The inputs of the filter are the meal and bolus (meal bolus and the extended bolus) information and the CGM data from a specific time in thje past until the current time.}
\label{fig:fig0}
\end{center}
\end{figure}

\begin{figure}[tb]
\begin{center}
\includegraphics[width=0.98\linewidth]{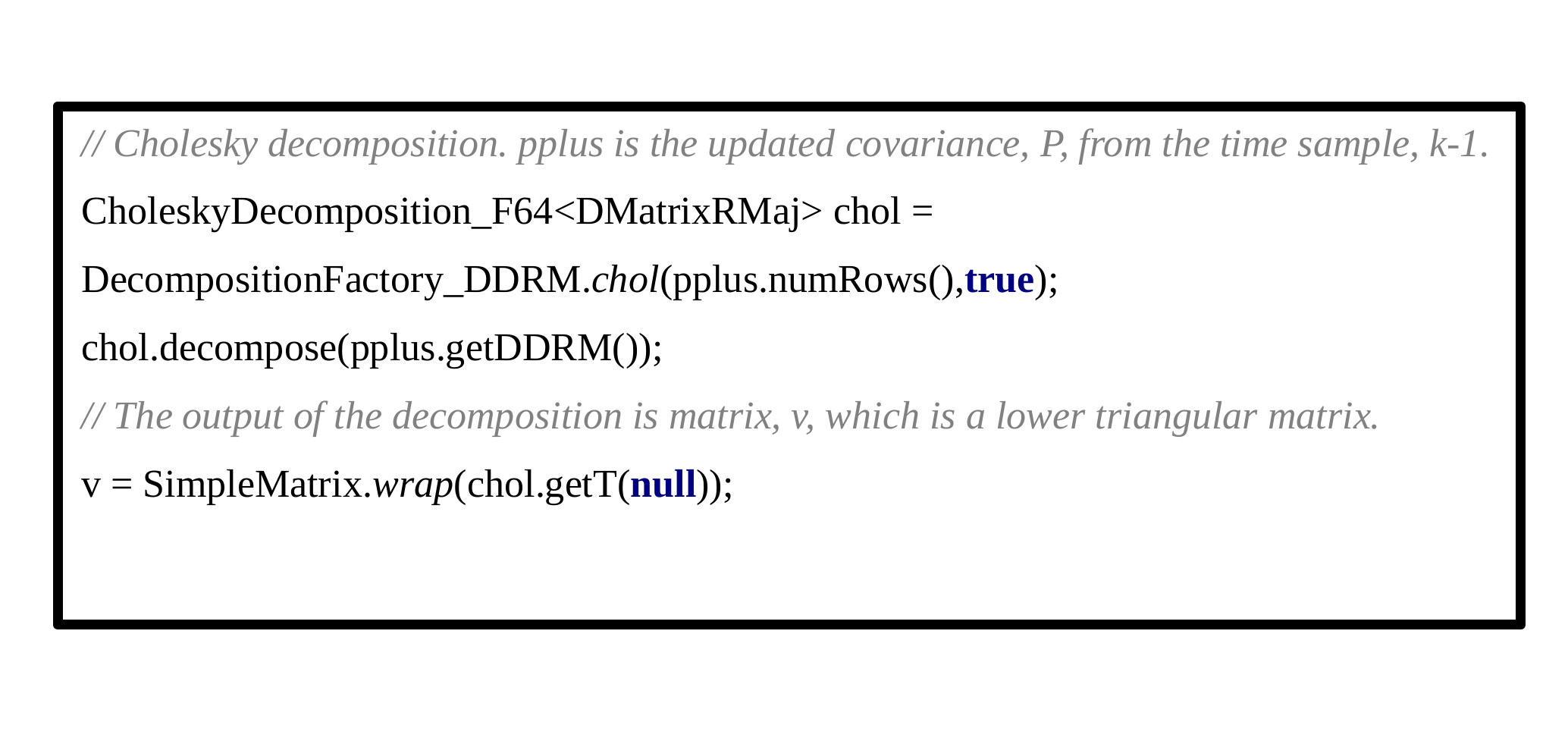}
\vspace{-0.5cm}
\caption{Cholesky decomposition for the UKF in EJML Java library.}
\label{fig:chol}
\end{center}
\end{figure}

Fig. \ref{fig:fig0} shows a few lines of the Java code from the implementation of the filters.
The basal insulin rate is directly added during the one-step prediction. The filtered data is then sent to the plot function. Fig. \ref{fig:chol} shows the code snippet for the Cholesky decomposition for the UKF in EJML.

\begin{figure}[tb]
\begin{center}
\includegraphics[width=0.98\linewidth]{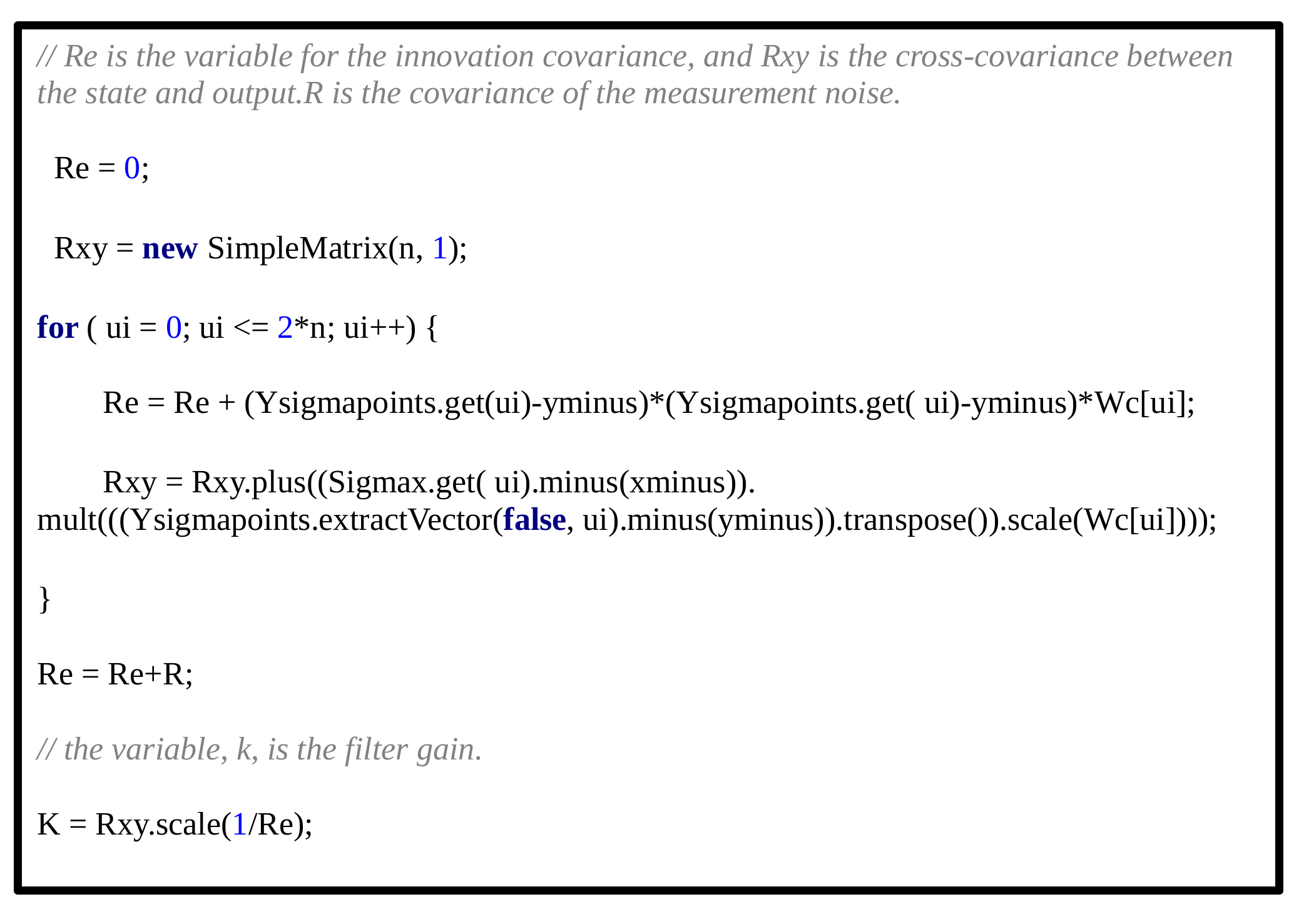}
\caption{Computation of innovation covariance, the cross-covariance in \eqref{rxy} and the filter gain in \eqref{k}, for the UKF in EJML Java library.}
\label{fig:rxy}
\end{center}
\end{figure}

Fig. \ref{fig:rxy} shows the code snippet for computation of the innovation covariance, the cross covariance, and the filter gain for the UKF in EJML.

\subsection{Hardware and data acquisition }

\begin{figure}[tb]
\begin{center}
\includegraphics[width=0.8\linewidth]{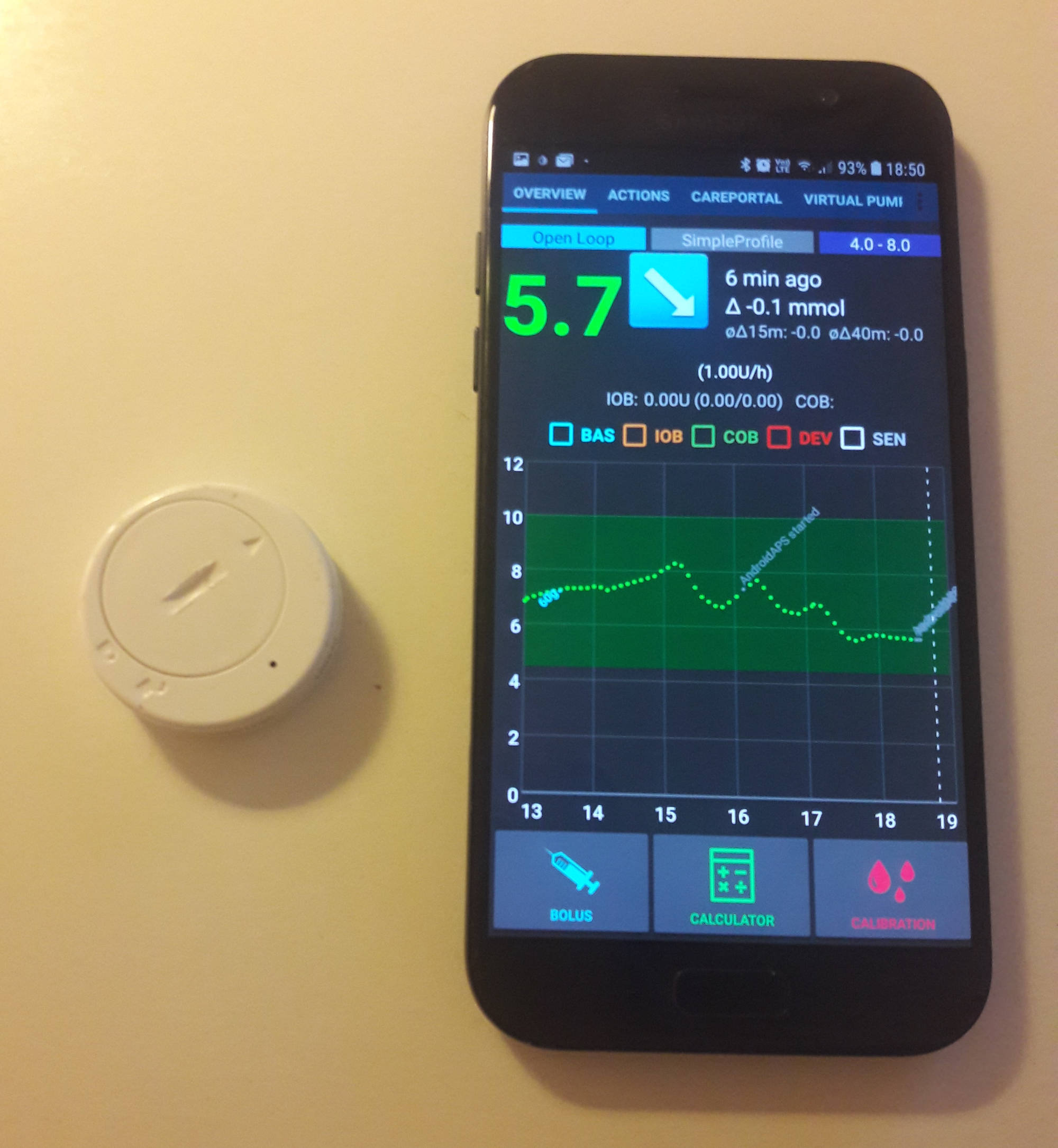}
\caption{Hardware used in the study consisting of FreeStyle Libre sensor, BluCon transmitter, and an Android smart phone. The BluCon transmitter is attached to the FreeStyle Libre CGM sensor. }
\label{fig:fig1}
\end{center}
\end{figure}

\begin{figure}[tb]
\vspace{0.5cm}
\begin{center}
\includegraphics[width=0.52\linewidth]{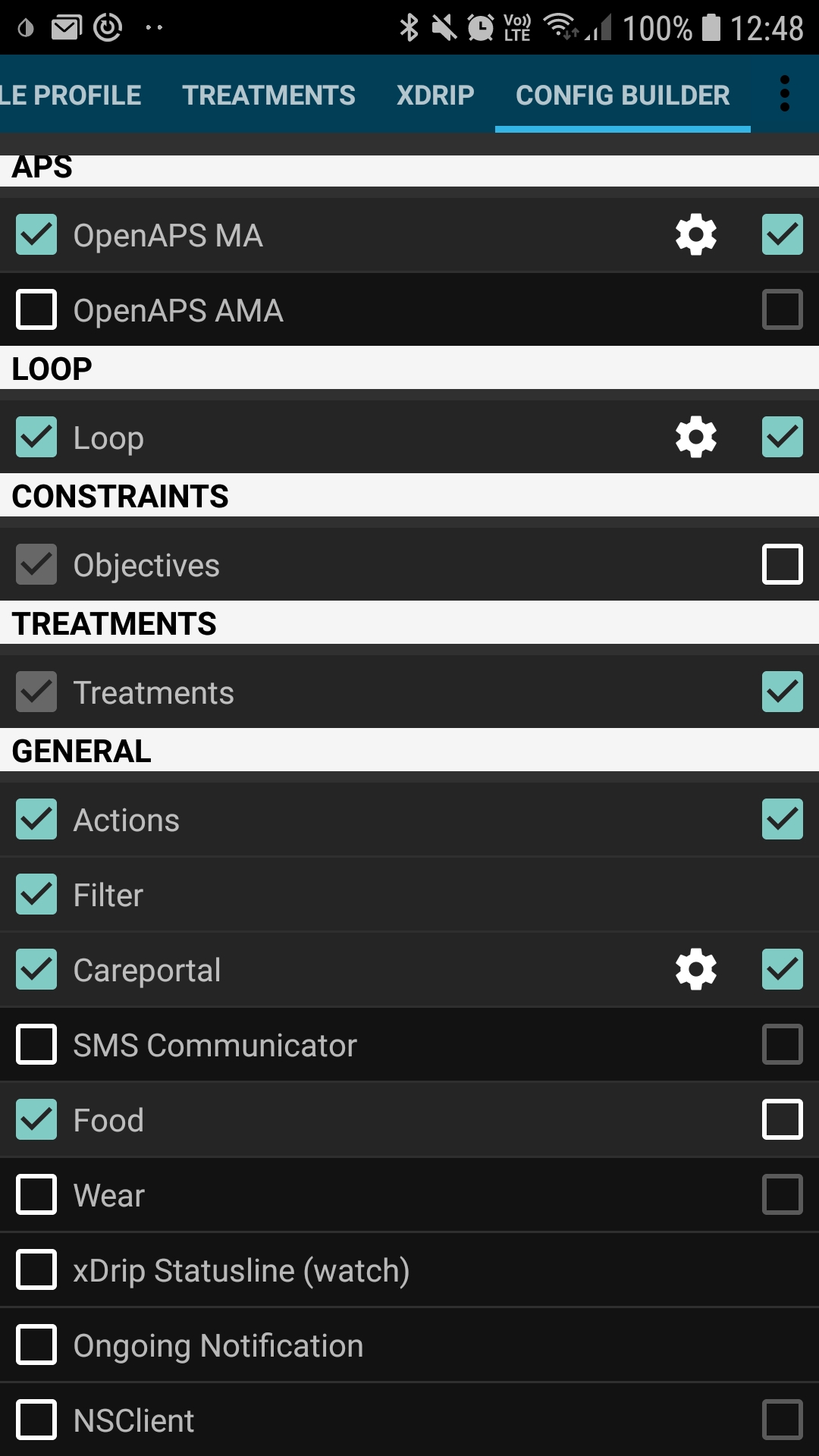}
\caption{The Filter implemented as a general Plugin in AndroidAPS.}
\label{fig:fig2}
\end{center}
\end{figure}

The filters are implemented in a Samsung A5 2017 running Android 8.1.0. We test implementation in Java of the filters in connection to a real CGM. To monitor and filter CGM data, we use AndroidAPS 1.58 (\url{http://github.com/MilosKozak/AndroidAPS}). AndroidAPS is an open source Java-based artificial pancreas application developed for Android phones. We added the filter as a plugin to AndroidAPS. The CGM is a FreeStyle Libre from Abbott Laboratories connected to a Blucon Bluetooth transmitter from Ambrosia Systems Inc. Fig. \ref{fig:fig1} illustrates the hardware configuration. Fig. \ref{fig:fig2} shows the filter implemented as a general plugin in AndroidAPS. 

\section{Results and discussion}
\label{sec:Results}

\subsection{Comparison of Java linear algebra libraries}

\begin{figure*}[!tb]
\centering
\subfigure[Matrix-matrix multiplication]{\includegraphics[width=0.49\linewidth]{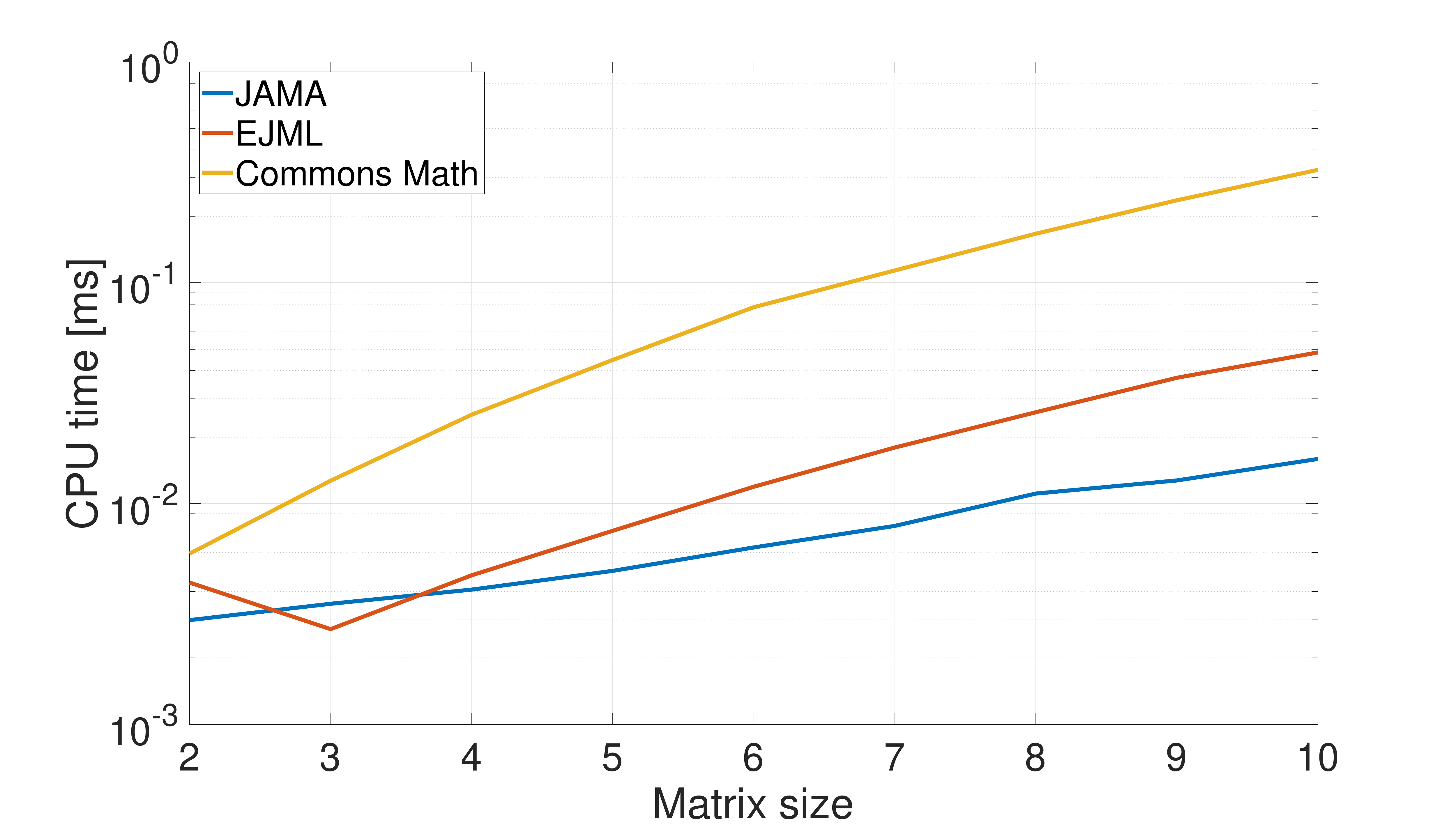} \label{fig:mult}}
\subfigure[Cholesky decomposition]{\includegraphics[width=0.49\linewidth]{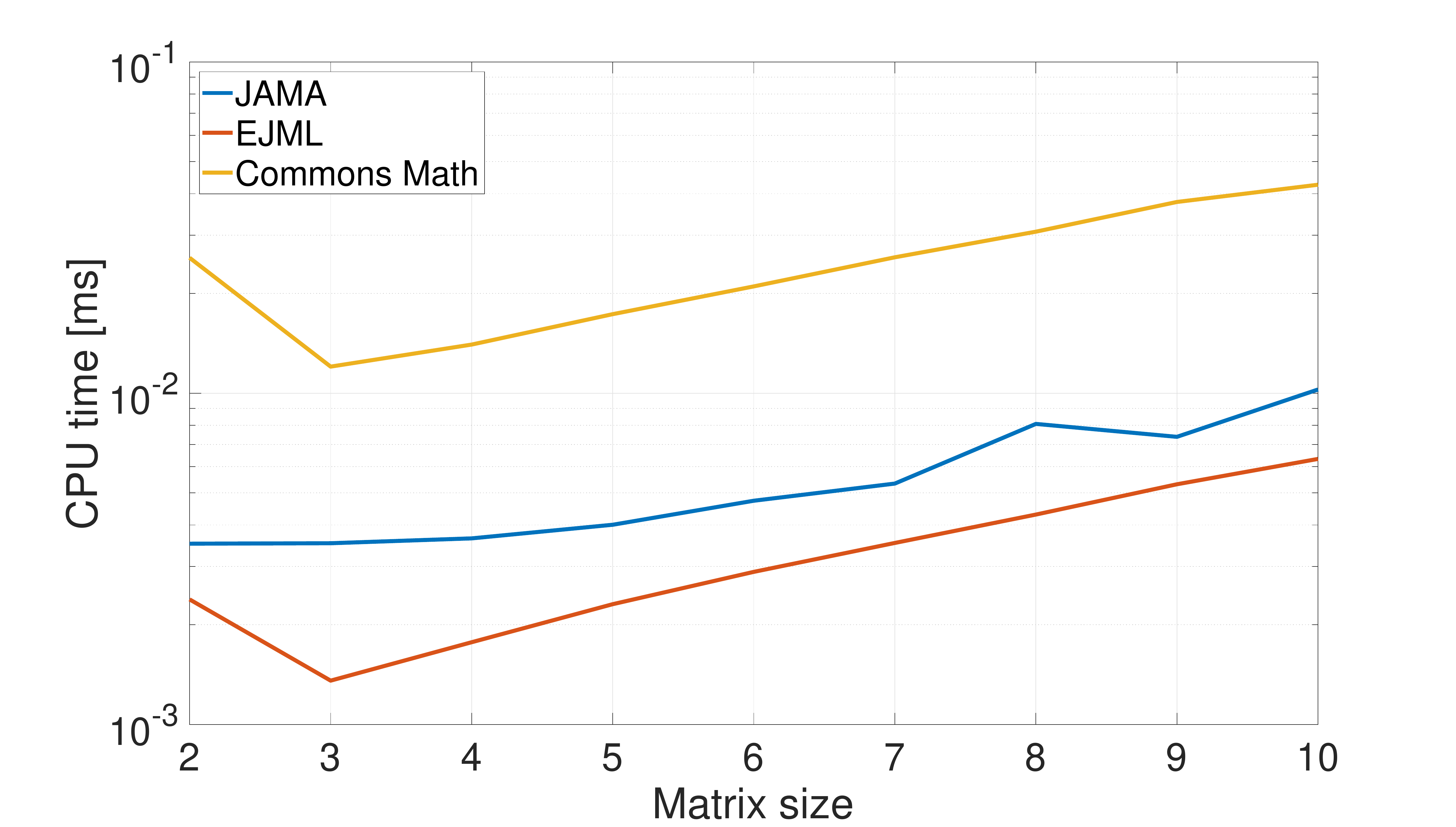}}
\subfigure[Matrix scaling]{\includegraphics[width=0.49\linewidth]{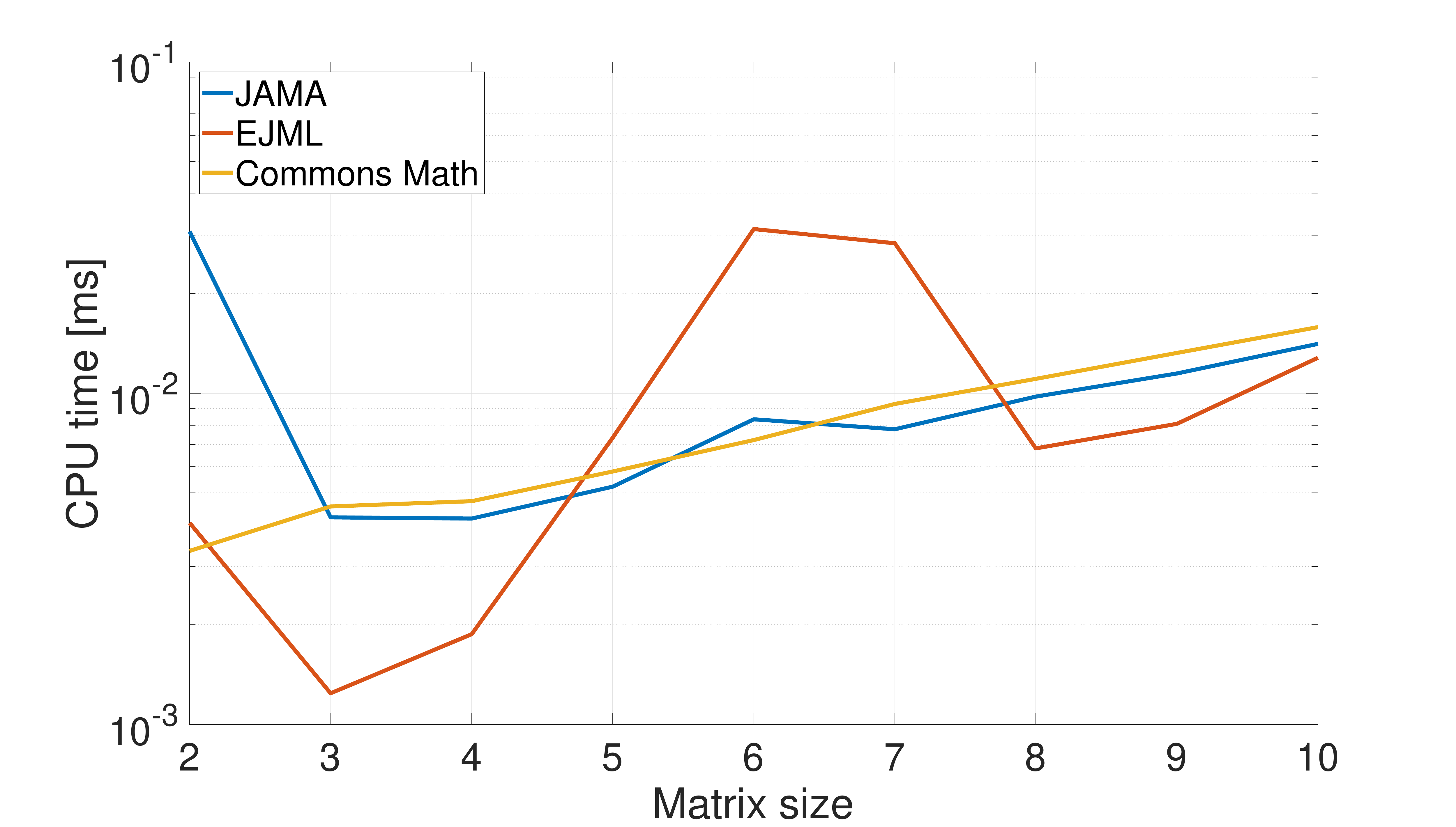}}
\subfigure[Matrix-matrix addition]{\includegraphics[width=0.49\linewidth]{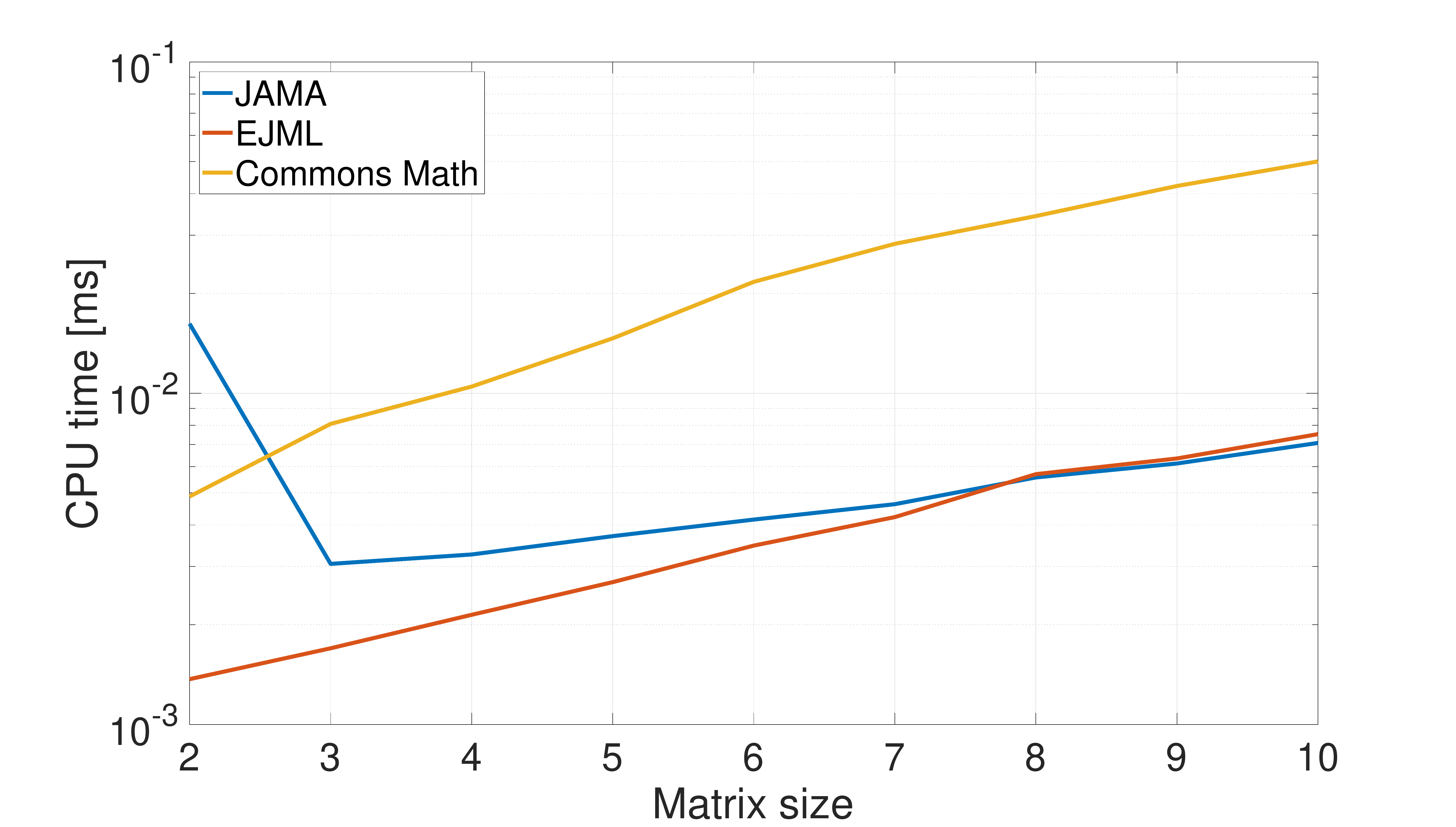}}
\caption{Comparison of computation time for five matrix operations in JAMA, EJML, and Commons math. The vertical axis is in logarithmic scale. The values are the average over 10000 runs.}
\label{fig:fig3}
\end{figure*}

Fig. \ref{fig:fig3} shows the computation time for four matrix algebraic operations in JAMA, EJML, and Commons math libraries in Java. The comparison indicates that EJML is faster than JAMA and Commons math libraries for Cholesky decomposition, and matrix-matrix addition. EJML is relatively comparable with the other two libraries for the matrix-matrix multiplication and matrix scaling. Therefore we choose to implement the filters in EJML. 

\subsection{Performance of linear and nonlinear Kalman filters}

\begin{figure}[!tb]
\begin{center}
\includegraphics[width=0.99\linewidth]{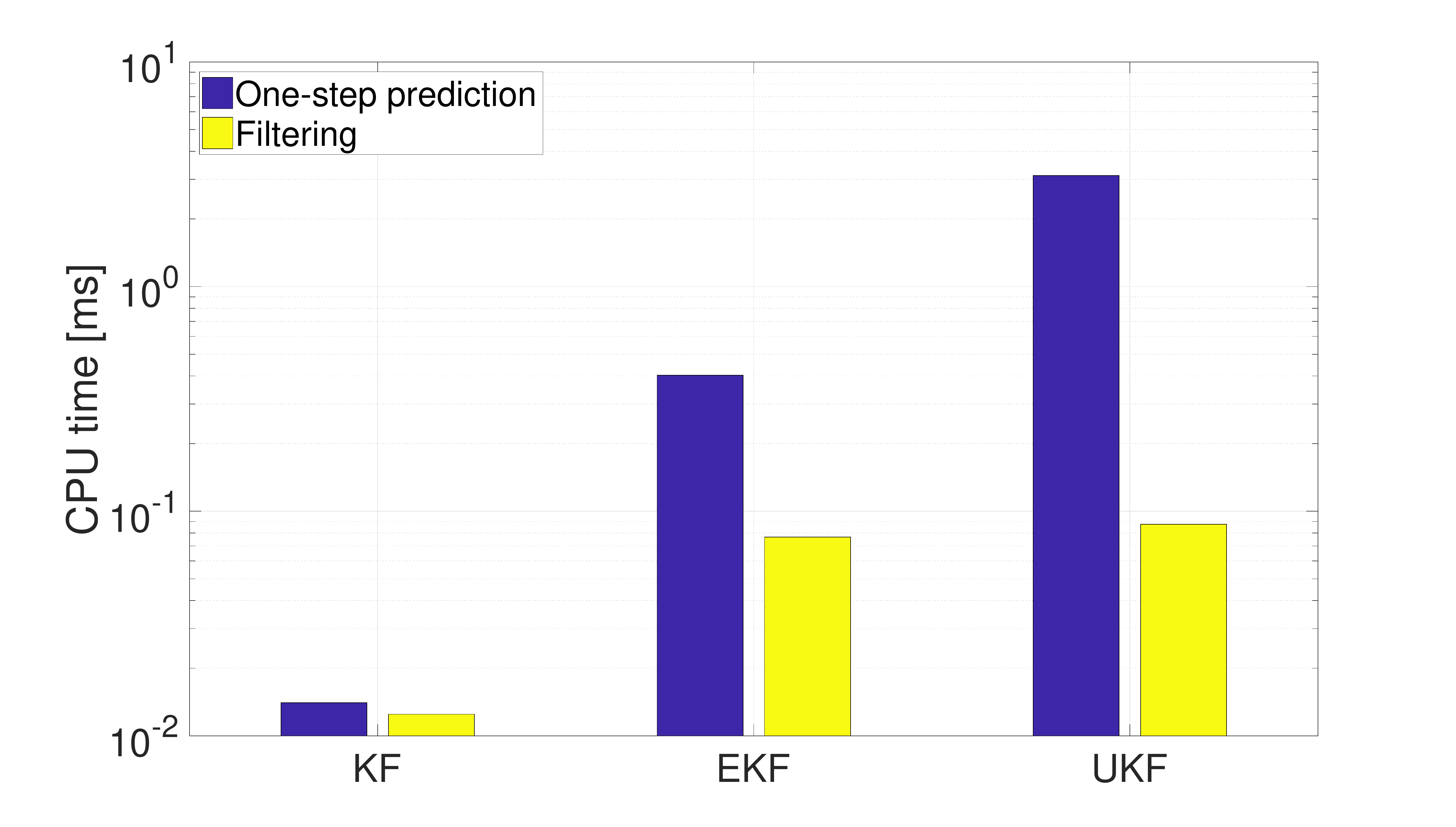}
\caption{Comparison of computation time for prediction step and filtering step, for Kalman filter, extended Kalman filter, and unscented Kalman filter. The values are the average over 10000 runs.}
\label{fig:fig4}
\end{center}
\end{figure}

Fig. \ref{fig:fig4} compares the processing time of the three filters in EJML. For the linear KF, we use a stationary filter, i.e. the state covariance matrix, the Kalman gain and the covariance of the innovation can be computed off-line. We use a forward Euler method to compute the one-step prediction of the state estimate for the UKF, and to compute the one-step prediction of the state covariance for the EKF. The step length is 1 minute, and the sampling time is 5 minutes, i.e. the one-step prediction requires five Euler steps.

Fig. \ref{fig:fig4} indicates that the linear KF has the smallest computation time for the prediction. This is expected as the one-step prediction in \eqref{predkf} does not need an Euler method implementation. In addition, the filtering step for the linear KF has the smallest computation time among the three filters, because it does not include a covariance update as the stationary state covariance is computed off-line outside the filtering step. The filtering computation time is similar for the EKF and UKF, while the one-step prediction for the UKF is more time-consuming than that for the EKF. For the UKF, it is required to compute the Cholesky factorization of $P_{k-1|k-1}$, the one-step predictions for the $2n+1$ sigma points and to compute the state covariance matrix, whereas the EKF requires to compute the one-step prediction of $n_x$ states, and the state covariance of size $n \times n$.

\begin{figure}[!tb]%
    \centering
    \subfigure[CGM filtered by the UKF]{{\includegraphics[width=4.3cm]{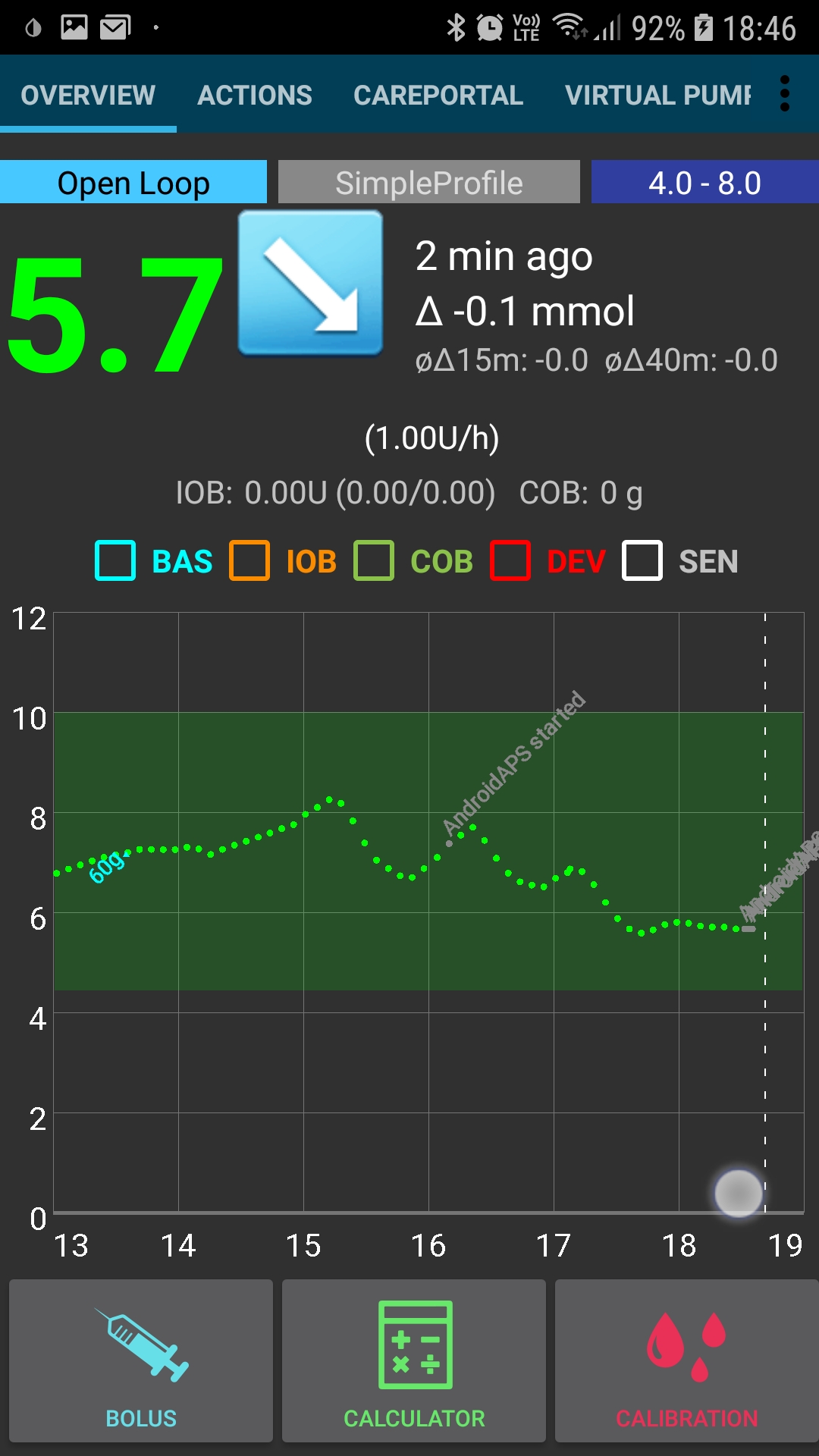} }}%
    \subfigure[Unfiltered CGM]{{\includegraphics[width=4.3cm]{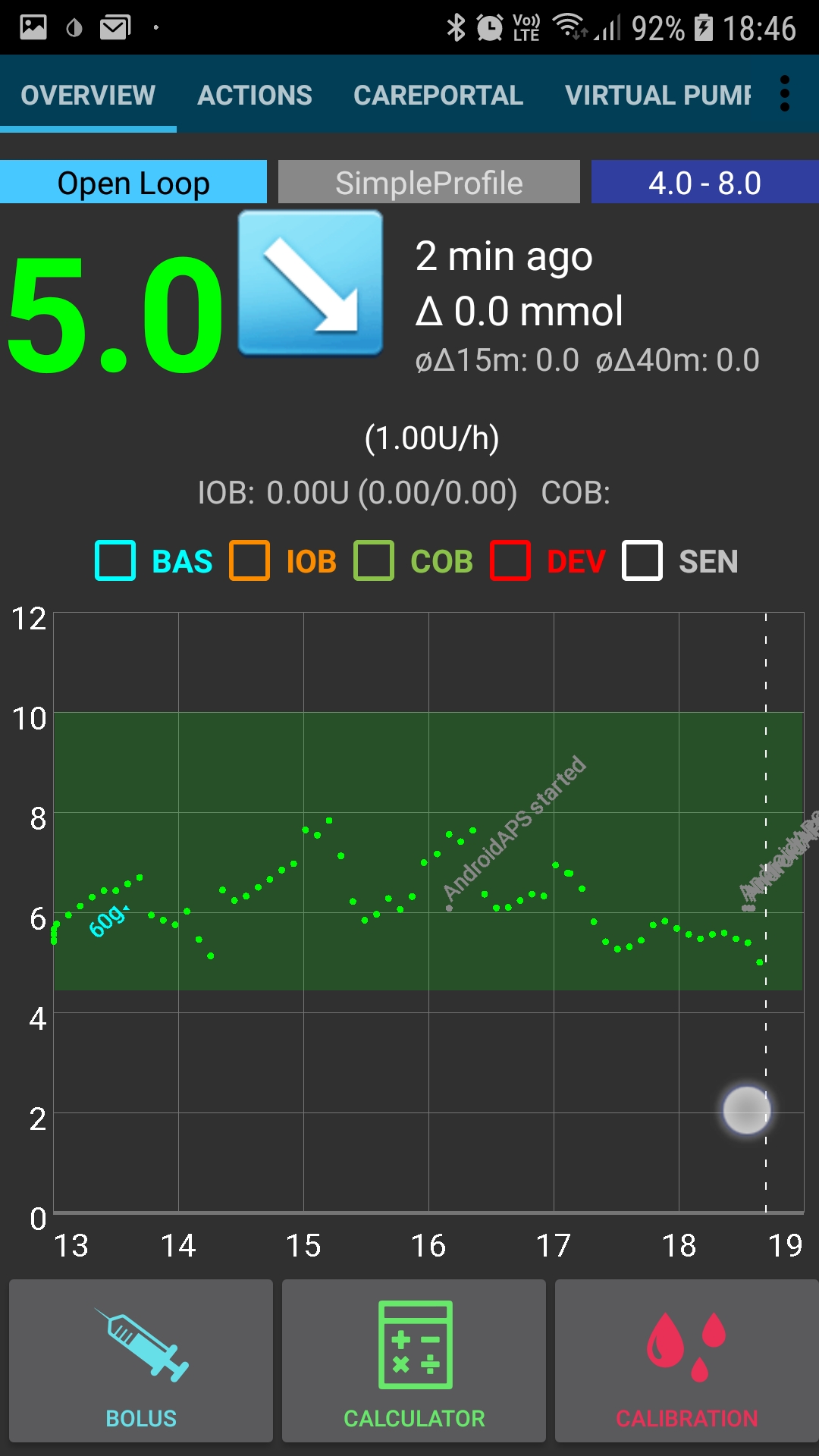} }}%
    \caption{The implemented UKF with the EJML library in AndroidAPS.}%
     \label{fig:fig5}%
\end{figure}

Fig. \ref{fig:fig5} shows the filtered and unfiltered CGM data. The filter is the UKF implemented in the Android smart phone using the EJML Java Matrix library in AndroidAPS platform. The filtered data is shown in the AndroidAPS graph when the Filter plugin in AndroidAPS is enabled (Fig. \ref{fig:fig2}). The unfiltered data are shown in the graph when the Filter plugin in AndroidAPS is disabled. The sensor is connected to a non-diabetes subject, while the model in the UKF is a diabetes patient model. Obviously, the model does not fit the CGM data and that is the reason for the relatively large bias between the unfiltered CGM data and the filtered data. This is of less importance in this study, because the aim of this implementation is to measure the computation time for the filters and solve software/hardware enginering issues with the Java implementation and the sensor connection.


\section{Conclusions}
\label{sec:Conclusions}
We presented the linear and nonlinear KF and implemented them in an Android smart phone using the Java programming language. We tested these filters in the AndroidAPS platform and use a FreeStyle Libre sensor as the source of CGM data. We compare the computation times of elementary matrix operations for three numerical linear algebra libraries implemented in Java. These operations are required for linear and nonlinear filtering. The results indicate that the EJML Java library is a suitable package for this purpose. The state-estimator in this study can be used for the implementation in Java of linear and nonlinear parameter and state estimation,  model-based controller, model-based sensor fault detection, as well as for meal detection and estimation. These filters and in particular the UKF can be used as part of an artificial pancreas, in model-based bolus advisors for connected insulin pens, and in meal-detection algorithms.


\end{document}